# Symmetric cube $L$-functions for $\mathrm{GL}_2$ are entire

By Henry H. Kim and Freydoon Shahidi *

## Introduction

The purpose of this paper is to prove the long awaited holomorphy of the third symmetric power $L$-functions attached to nonmonomial cusp forms of $\mathrm{GL}_2$ over an arbitrary number field on the whole complex plane.

To be more precise, let $F$ be a number field whose ring of adeles is $\mathbb{A} = \mathbb{A}_F$. Let $\pi = \otimes_v \pi_v$ be a cuspidal (unitary) representation of $\mathrm{GL}_2(\mathbb{A})$. Let $S$ be a finite set of places of $F$ such that for $v \notin S$, $\pi_v$ is unramified. For each $v \notin S$, let

$$t_v = \left\{ \begin{pmatrix} \alpha_v & 0 \\ 0 & \beta_v \end{pmatrix} \right\}$$

denote the semisimple conjugacy class of $\mathrm{GL}_2(\mathbb{C})$ defining $\pi_v$ (cf. [Bo], [La4]). We recall that if $\pi$ is attached to a classical modular form of weight $k$ on the upper half-plane for which the Fourier coefficient at $p$ is $a_p$, then

$$\begin{aligned} a_p &= p^{\frac{k-1}{2}}(\alpha_p + \beta_p) \\ &= p^{\frac{k-1}{2}}(\alpha_p + \alpha_p^{-1}). \end{aligned}$$

Fix a positive integer $m$. Denote by $r_m = \mathrm{sym}^m(\rho_2)$, the $m^{\mathrm{th}}$ symmetric power representation of the standard representation $\rho_2$ of $\mathrm{GL}_2(\mathbb{C}) = {}^L\mathrm{GL}_2^0$, an irreducible representation of dimension $m+1$ (cf. [De], [La4], [Se], [Sh5, 7], [T]). Set

$$L_S(s, \pi, r_m) = \prod_{v \notin S} L(s, \pi_v, r_m),$$

where

$$\begin{aligned} L(s, \pi_v, r_m) &= \det(I - r_m(t_v)q_v^{-s})^{-1} \\ &= \prod_{0 \leq j \leq m} (1 - \alpha_v^j \beta_v^{m-j} q_v^{-s})^{-1} \end{aligned}$$

*The first author was partially supported by NSF grant DMS9610387. The second author was partially supported by NSF grant DMS9622585.

is the Langlands $L$-function attached to $\pi_v$ and $r_m$. Recall that if $O_v$ and $P_v$ are the ring of integers of $F_v$ and its unique maximal ideal, then $q_v = \text{card}(O_v/P_v)$.

The analytic properties (cf. [Sh7]) and special values (cf. [De]) of these $L$-functions are of great interest and their consequences have been noted by several mathematicians [De], [La4], [Muk], [Mur], [Se], [Sh7], [T].

A cuspidal representation $\pi$ of $\text{GL}_2(\mathbb{A})$ is called monomial if $\pi \cong \pi \otimes \eta$ for some nontrivial grossencharacter $\eta$, necessarily quadratic. Let $K/F$ be the quadratic extension of $F$ attached to $\eta$ by class field theory. Then there exists a quasicharacter $\chi = \otimes_w \chi_w$ of $K^* \backslash \mathbb{A}_K^*$ which does not factor through the norm so that if $\sigma_v = \text{Ind}(W_{F_v}, W_{K_w}, \chi_w)$, then $\pi_v = \pi(\sigma_v)$, the representation attached to $\sigma_v$ by Langlands correspondence or the Weil representation attached to $\chi_w$. Here $W_*$ denotes the corresponding Weil group, $w|v$, and $\pi = \otimes_v \pi_v$ (cf. [L-La], [S-Tan]).

Throughout this paper we shall assume $m = 3$. By Langlands [La1], $L_S(s, \pi, r_3)$ is a meromorphic function of $s$ on $\mathbb{C}$. If $\pi \cong \pi \otimes \eta$, $\eta \neq 1$, then the completed $L$-function

$$L(s, \pi, r_3) = L(s, \chi^3) L(s, \chi^2 \chi'),$$

where the $L$-functions on the right are those of Hecke attached to quasicharacters of $K^* \backslash \mathbb{A}_K^*$ (cf. Proposition 2.3). The quasicharacter $\chi'$ is the conjugate of $\chi$ by the nontrivial element of the Galois group. Without loss of generality we may assume the central character $\omega_\pi = \eta \cdot \chi|\mathbb{A}^*$ is trivial on $\mathbb{R}_+^*$, where $\mathbb{A}^* = \mathbb{I} = \mathbb{I}^1 \cdot \mathbb{R}_+^*$ with $\mathbb{I}^1$ ideles of norm 1 and $\mathbb{R}_+^*$ the multiplicative group of positive real numbers. Observe that this is equivalent to $\chi|\mathbb{R}_+^* \equiv 1$. Thus $L(s, \pi, r_3)$ for such $\pi$ has poles only at $s = 0$ and 1 if and only if $\chi^3 = 1$. Moreover, the poles are simple. Observe that since $\pi$ is cuspidal $\chi \neq 1$. In this paper we prove (Theorem 4.9):

THEOREM 1. *Suppose $\pi$ is not monomial. Then the partial symmetric cube $L$-function $L_S(s, \pi, r_3)$ is entire.*

In fact, what we prove is a much stronger result: the holomorphy of the full (completed) $L$-function on all of $\mathbb{C}$, stated in this introduction as Theorem 3, from which Theorem 1 follows.

The interest in this $L$-function and such a result has a long history and has attracted the attention of number theorists (cf. [D-I], [De], [H-Ra], [Sh5, 7], for example). Ever since Shimura [Shi] and Gelbart-Jacquet [Ge-J] established similar results for $L_S(s, \pi, r_2)$ which led to the so-called Gelbart-Jacquet or adjoint square lift of $\pi$ to $\text{GL}_3(\mathbb{A})$ with numerous important applications (cf. [B-D-H-I], [H-Ra], [Lu-R-Sa1], [Lu-R-Sa2], [La3], [Ra], [Sh4, 7], to name a few), similar questions were posed for $L_S(s, \pi, r_3)$ with the hope that, beside its arithmetic consequences, it could lead to a lift of $\pi$ to $\text{GL}_4(\mathbb{A})$, a result of great importance (cf. [Sh7] and in particular the estimate 3/22 for Fourier coefficients of Maass forms).

It was Langlands who first established the meromorphy of $L_S(s, \pi, r_3)$ on $\mathbb{C}$ by expressing it as the constant term of an Eisenstein series on an exceptional group of type $G_2$ as one of many examples in his *Euler Product* monograph [La1], [Go]. His ideas about their functional equations which were explained in a letter to Godement, were taken up by Shahidi [Sh2, 3, 1], [C-S] who established their functional equations in general. This led to a new approach to the study of automorphic L-functions which is now referred to as the Langlands-Shahidi method [La1], [M-W2], [Sh1, 3, 4, 5], [K3], [F-Gol]. In the case of $L_S(s, \pi, r_3)$, similar results were observed by Deligne who used the Gelbart-Jacquet lift $\Pi$ of $\pi$ together with the properties of $L(s, \pi \times \Pi)$ (cf. [J-S], [M-W2], [Sh5], and Proposition 2.2).

The first instances where the holomorphy of $L_S(s, \pi, r_3)$ on all of $\mathbb{C}$ were established were proved in [Mo-Sh] and [Sh5], and it was in [Sh5] that the local factors $\varepsilon(s, \Pi_v, r_3, \psi_v)$ and $L(s, \pi_v, r_3)$ were defined for all $v$ and the full functional equation was stated. Moreover the adjoint cube representation $r_3^0 = r_3 \otimes (\Lambda^2 r_1)^{-1}$ was defined and the machinery for applying the converse theorem using $L(s, \pi, r_3^0, \rho)$, $\rho$ a grossencharacter, was set up (Theorem 4.1 of [Sh5]). But the results were conditional (cf. Proposition 2.2.2 here and Theorem 4.1 of [Sh5]) and one purpose of the present paper is to establish the result with no restriction (Theorem 2 below).

Although several efforts were made, leading to new and interesting L-functions, it was not until 1994 (date of the preprint) that the first integral representation for $L_S(s, \pi, r_3)$ was found by Bump, Ginzburg, and Hoffstein [B-G-H]. Their integral generalized that of Shimura [Shi] and Gelbart-Jacquet [Ge-J] and therefore required a 3-sheet cover of $\mathrm{GL}_2$, necessitating that $F$ contain 3rd roots of 1. It was there that holomorphy of $L_S(s, \pi, r_3)$ up to $\mathrm{Re}(s) > 3/4$ was established, extending the bound $\mathrm{Re}(s) > 1$ of absolute convergence in [Sh5], a consequence of a similar statement for $L(s, \pi \times \Pi)$ due to Jacquet and Shalika [J-S]. Uncompromising complications at archimedean places have stopped the extension beyond 3/4.

Our method is that of Langlands-Shahidi and is therefore quite representation theoretic. It starts with an observation of Kim (Observation 4.1), used first in [K3] to prove that the exterior square L-functions for $\mathrm{GL}_n$, $n$ odd, are entire, and a consequence of Langlands' theory of Eisenstein series [La2], [M-W1], that if $L(s, \pi, r_3^0)$ has a pole for $1/2 < s < 1$, then the corresponding residual representation must be unitary [K2]. We recall that arguments of this type, but in the opposite direction, were used by Speh [Sp] and later by Tadic [Ta] to determine unitary duals of local groups. Since the unitary dual of $G_2$ over a $p$-adic field (for real groups see [V]) is now completely determined (see a recent article of Muić [Mu]), this gives a contradiction, but not until we bring an indispensable global element into our proof. (See the discussion before Theorem 4.3.) This is provided by Ramakrishnan's recent result [Ra] which states that many of the $\pi_v$, $v \notin S$, are tempered (true for all $\pi_v$ by Deligne, if $\pi$ is holomorphic). We remark that the possibility of poles on $\text{Re}(s) = 1/2$ needs to be eliminated by other, all global, means (cf. [Sh5], [La2], and Proposition 2.2 and Lemma 4.7 here).

We give several proofs for Theorem 1 (Theorem 4.9), the first two relying in varying degrees on [Ra]. The first relies more on global information [Ra], while the second, which uses the hermitian property of each $J(s, \pi_v)$ if $s$ is a pole (Observation 4.4), relies more on the classification of Muić [Mu]. As expected, neither proof is completely local and in fact the first proof uses the existence of tempered components at a set of (lower) density $> 1/2$. The possibility of a pole at $s = 1/2$ in the first proof is eliminated by using a quadratic base change [La3], the holomorphy at $s = 1/2$ if the central character is trivial [Sh5], a result also used to eliminate the pole at $s = 1/2$ in the second proof, and the functional equation [Sh2, 5]. The theorem now follows from the functional equation. (See Section 2.)

We should remark that having even one unramified tempered component $\pi_v$ for an arbitrary cuspidal representation $\pi$ of $\text{GL}_2(\mathbb{A})$ which was used in the second proof, is highly nontrivial. In fact, if $n \geq 3$, although expected, nothing as such is known, even though it appears that the methods of [Ra] will give the existence of tempered components for self contragredient forms on $\text{GL}(3)$.

Third, we sketch a proof in which our global ingredient is much weaker than the two others, namely: Given a cusp form $\pi$ and any $\varepsilon > 0$, there exists a local component $\pi_v = \pi(\mu_v | \ |_v^{r_v}, \mu_v | \ |_v^{-r_v})$, $0 \leq r_v < 1/2$, $\mu_v$ a character of $F_v^*$, such that $r_v < \varepsilon$ (cf. [K3]). But, we then need to make much deeper use of the unitary dual of $G_2$ coming from the conjugacy class of Borel subgroups ([Mu, Th. 5.2]). This is the approach pursued in [K3] as one does not expect strong results such as those in [Ra] to be available when more general cases are treated.

As explained in [Sh5], to write the result in a language appropriate for applying the converse theorem [Co-PS], one needs to consider $r_3^0 = r_3 \otimes (\Lambda^2 r_1)^{-1}$. Let $\rho = \otimes_v \rho_v$ be a grossencharacter. If $L(s, \pi_v, r_3^0, \rho_v)$ is the corresponding local $L$-function for each $v$ as defined in [Sh5], i.e. simply

$$L(s, \pi_v, r_3^0, \rho_v) = L(s, \pi_v \otimes \rho_v, r_3^0),$$

let

$$L(s, \pi, r_3^0, \rho) = \prod_v L(s, \pi_v, r_3^0, \rho_v).$$

Similarly define a root number $\varepsilon(s, \pi, r_3^0, \rho)$ as in [Sh5]. We have:

THEOREM 2. a) *The L-function $L(s, \pi, r_3^0, \rho)$ is entire unless $\pi$ is monomial. Suppose $\pi \cong \pi \otimes \eta$, $\eta \neq 1$. Let $K/F$ be the quadratic extension attached to $\eta$. Assume $\pi = \pi(\sigma)$, where $\sigma = \mathrm{Ind}(W_F, W_K, \chi)$. Then*

$$L(s, \pi, r_3^0, \rho) = L(s, \chi^2 {\chi'}^{-1} \rho_K) L(s, \chi \rho_K),$$

$\rho_K = \rho \cdot N_{K/F}$, *and if $\omega_\pi \rho^2 | \mathbb{R}_+^* \equiv 1$, then the L-function $L(s, \pi, r_3^0, \rho)$ has poles only at $s = 0$ and $1$. They appear if and only if $\chi^2 {\chi'}^{-1} \rho_K = 1$ and are simple.*
b) *The standard functional equation*

$$L(s, \pi, r_3^0, \rho) = \varepsilon(s, \pi, r_3^0, \rho) L(1-s, \tilde{\pi}, r_3^0, \rho^{-1})$$

*is satisfied.*

One still needs to twist with arbitrary cusp forms on $\mathrm{GL}_2(\mathbb{A})$ (cf. Theorem 8.2 of [Sh7]), but with this theorem we have hopefully taken a major step in establishing the symmetric cube lift from $\mathrm{GL}_2(\mathbb{A})$ to $\mathrm{GL}_4(\mathbb{A})$, i.e. the existence of an automorphic form on $\mathrm{GL}_4(\mathbb{A})$ whose standard L-function is equal to $L(s, \pi, r_3)$.

More precisely, if $\Pi$ is the adjoint cube lift of $\pi$, which as explained one expects to establish using a converse theorem, then the symmetric cube lift of $\pi$ is simply $\Pi \otimes \omega_\pi$, where $\omega_\pi$ is the central character of $\pi$ (cf. [Sh5]).

We conclude by stating the full result for $r_3$ from which Theorem 1 follows. By substituting $\pi$ with $\pi \otimes \omega_\pi$ and setting $\rho = 1$ in Theorem 2 we have (Theorem 4.9):

THEOREM 3. a) *The L-function $L(s, \pi, r_3)$ is entire unless $\pi$ is monomial. Suppose $\pi \cong \pi \otimes \eta$, $\eta \neq 1$. Let $K/F, \chi$, and $\sigma$ be as in Theorem 2 so that $\pi = \pi(\sigma)$. Then*

$$L(s, \pi, r_3) = L(s, \chi^3) L(s, \chi^2 \chi')$$

*and the poles are simple and if $\omega_\pi^3 | \mathbb{R}_+^* \equiv 1$, they appear only at $s = 0$ and $1$. The poles exist if and only if $\chi^3 = 1$.*
b) *The standard functional equation*

$$L(s, \pi, r_3) = \varepsilon(s, \pi, r_3) L(1-s, \tilde{\pi}, r_3)$$

*is satisfied.*

As a consequence of Theorem 3, there is another proof of the following result. We refer to [Ga], [Ik], [PS-Ra1] for original proofs (cf. Remarks 1.1 and 4.12).

COROLLARY. *Suppose $\pi$ is not monomial. Then the partial Rankin triple L-function $L_S(s, \pi \times \pi \times \pi)$ is entire.*

We should point out that Theorem 3 can also be used to supplement the necessary holomorphy condition needed in the main Theorem C of [H-Ra]

which was provided by [B-G-H] there. One can then dispose of Lemma 5.4 of [H-Ra].

The authors hope to take up related arithmetic questions in future papers.

We would like to thank Dinakar Ramakrishnan for his careful reading of the manuscript and many useful comments.

## 1. Preliminaries

Let $F$ be a number field and let $G$ be a split group of type $G_2$ defined over $F$. For each place $v$ of $F$ let $F_v$ be the corresponding completion and if $v$ is finite, let $O_v$ be its ring of integers. If $P_v$ is the maximal ideal of $O_v$, set $q_v$ for the cardinability of $O_v/P_v$. If $\mathbb{A}$ is the ring of adeles of $F$, we set $\mathbb{A}_\infty = \prod_{v=\infty} F_v$. Fix a split Cartan subgroup $T$ in $G$ and let $B = TU$ be a Borel subgroup of $G$. Let $K_\infty$ be the standard maximal compact subgroup in $G(\mathbb{A}_\infty)$ and let $K_v = G(O_v)$ for a finite place $v$. The product $K = K_\infty \times \prod_{v<\infty} K_v$ is a maximal compact subgroup in $G(\mathbb{A})$.

We follow Moeglin and Waldspurger [M-W1, App. 2]. In what follows the roots are those of $T$ in $U$. Let $\beta_1$ be the long simple root and $\beta_6$ the short one. Let

$$\beta_2 = \beta_1 + \beta_6, \quad \beta_3 = 2\beta_1 + 3\beta_6, \quad \beta_4 = \beta_1 + 2\beta_6, \quad \beta_5 = \beta_1 + 3\beta_6,$$

denote the other positive roots.

Let $P$ the maximal parabolic subgroup generated by $\beta_1$ (the long root). Then we have a Levi decomposition ([Sh5]) $P = MN$, with $M \simeq \mathrm{GL}_2$. Thus

$$\mathfrak{a}^* = X(M) \otimes \mathbb{R} = \mathbb{R}\beta_4, \text{ and } \mathfrak{a} = \mathbb{R}\beta_4^\vee.$$

Here $X(M)$ is the group of $F$-rational characters of $M$. If $\rho_P$ is half the sum of roots generating $N$, then $\rho_P = \frac{5}{2}\beta_4$.

Let $\tilde{\alpha} = \beta_4$ and identify $s \in \mathbb{C}$ with $s\tilde{\alpha} \in \mathfrak{a}_\mathbb{C}^*$. Then $s\tilde{\alpha}$ corresponds to the character $|\det(m)|^s$. Let $\pi = \otimes \pi_v$ be a cusp form on $M(\mathbb{A}) = \mathrm{GL}_2(\mathbb{A})$. We may and will assume that the central character $\omega_\pi$ of $\pi$ is trivial on $\mathbb{R}_+^*$, where $\mathbb{A}^* = \mathbb{I} = \mathbb{I}^1 \cdot \mathbb{R}_+^*$ with $\mathbb{I}^1$ ideles of norm 1. Given a $K$-finite function $\varphi$ in the space of $\pi$, we extend $\varphi$ to a function $\tilde{\varphi}$ on $G(\mathbb{A})$ as in [Sh3] and set

$$\Phi_s(g) = \tilde{\varphi}(g) \exp\langle s + \rho_P, H_P(g)\rangle.$$

Here $H_P$ is the usual Harish-Chandra homomorphism of $M$ into $\mathrm{Hom}(X(M), \mathbb{R})$. We define an Eisenstein series by

$$E(s, \tilde{\varphi}, g) = \sum_{\gamma \in P(F)\backslash G(F)} \Phi_s(\gamma g).$$

It is known that $E(s, \tilde{\varphi}, g)$ converges for $\mathrm{Re}(s) \gg 0$ and extends to a meromorphic function of $s$ on $\mathbb{C}$, with a finite number of poles in the plane $\mathrm{Re}(s) > 0$, all simple and on the real axis.

It is also known that $L^2_{\mathrm{dis}}(G(F)\backslash G(\mathbb{A}))_{(M,\pi)}$ is spanned by the residues of Eisenstein series for $\mathrm{Re}(s) > 0$. We note that for each $s$, the representation of $G(\mathbb{A})$ on the space of $\Phi_s$ is equivalent to $I(s,\pi) = \mathrm{Ind}_P^G \pi \otimes \exp(\langle s, H_P(\ )\rangle)$, where $H_P$ is the Harish-Chandra homomorphism. We know that the poles of Eisenstein series coincide with those of its constant terms. So it is enough to consider the constant term along $P$, which is given for each $f \in I(s,\pi)$ by

$$E_0(s,f,g) = \sum_{w \in \Omega} M(s,\pi,w)f(g),$$

$\Omega = \{1, \rho_6\rho_1\rho_6\rho_1\rho_6\}$, and

$$M(s,\pi,w)f(g) = \int_{N_w^-(\mathbb{A})} f(w^{-1}ng)\,dn.$$

Weyl group representatives are all chosen to lie in $K \cap G(F)$. Here

$$N_w^- = \prod_{\substack{\alpha > 0 \\ w^{-1}\alpha < 0}} U_\alpha,$$

where $U_\alpha$ is the one parameter unipotent subgroup. Then

$$M(s,\pi,w) = \otimes M(s,\pi_v,w),$$

where

$$M(s,\pi_v,w)f_v(g) = \int_{N_w^-(F_v)} f_v(w^{-1}ng)\,dn$$

with $f = \otimes f_v$, where $f_v$ is the unique $K_v$-fixed function normalized by $f_v(e_v) = 1$ for almost all $v$.

Let $^LM = \mathrm{GL}_2(\mathbb{C})$ be the $L$-group of $M$. Denote by $r$ the adjoint action of $^LM$ on the Lie algebra $^L\mathfrak{n}$ of $^LN$, the $L$-group of $N$.

Then

$$r = r_3^0 \oplus \wedge^2 r_1,$$

where $r_1$ is the standard representation of $\mathrm{GL}_2(\mathbb{C})$, $r_3$ is the symmetric cube representation of $\mathrm{GL}_2(\mathbb{C})$, and $r_3^0 = r_3 \otimes (\wedge^2 r_1)^{-1}$ is the adjoint cube representation of $\mathrm{GL}_2(\mathbb{C})$ (see [Sh5]).

Then it is well-known ([La1], [Sh4], equation (2.7)) that for $w = \rho_6\rho_1\rho_6\rho_1\rho_6$

$$M(s,\pi,w)f = \bigotimes_{v \in S} M(s,\pi_v,w)f_v \otimes \bigotimes_{v \notin S} \tilde{f}_v \times \frac{L_S(s,\tilde{\pi},r_3^0)L_S(2s,\tilde{\pi},\wedge^2 r_1)}{L_S(1+s,\tilde{\pi},r_3^0)L_S(1+2s,\tilde{\pi},\wedge^2 r_1)}$$

where $S$ is a finite set of places of $F$, including all the archimedean places such that for every $v \notin S$, $\pi_v$ is a class one representation and if $f = \otimes_v f_v$, then for $v \notin S$, $f_v$ is the unique $K_v$-fixed function normalized by $f_v(e_v) = 1$. The function $\tilde{f}_v$ is the $K_v$-fixed function in the space of $I(-s, w(\pi_v))$.

Finally, $L_S(s, \pi, r_3^0) = \prod_{v \notin S} L(s, \pi_v, r_3^0)$, where $L(s, \pi_v, r_3^0)$ is the local Langlands' $L$-function attached to $\pi_v$ and $r_3^0$, and $L_S(s, \pi, \wedge^2 r_1) = L_S(s, \omega_\pi)$ is the partial Hecke $L$-function, where $\omega_\pi$ is the central character of $\pi$. Finally $\tilde{\pi}$ is the contragredient of $\pi$ and $\omega_{\tilde{\pi}} = \omega_\pi^{-1}$.

*Remark* 1.1. Ikeda [Ik] calculated the poles of the Rankin triple $L$-function $L_S(s, \pi \times \pi \times \pi)$ for $\pi$ a cuspidal representation of $\mathrm{GL}_2$ [Ga], [PS-Ral]. The triple $L$-function is related to the partial symmetric cube $L$-function of $\pi$ as follows:

$$L_S(s, \pi \times \pi \times \pi) = L_S(s, \pi, r_3)(L_S(s, \pi \otimes \omega_\pi))^2.$$

The symmetric cube $L$-function satisfies

$$L(s, \pi, r_3) = L(s, \pi \otimes \omega_\pi, r_3^0).$$

## 2. A review of earlier global results

We will continue with our assumption that $\omega_\pi|\mathbb{R}_+^* \equiv 1$. As mentioned before we will lose no generality and therefore we resume this assumption throughout our proof. We know the following facts about the symmetric cube $L$-functions. Let $S$ be a finite set of places of $F$ such that every $\pi_v$, $v \notin S$, is of class one, i.e., it has a vector fixed by $\mathrm{GL}_2(O_v)$. Then the class of each $\pi_v$, $v \notin S$, is uniquely determined by the conjugacy class of a diagonal matrix $t_v = \begin{pmatrix} \alpha_v & 0 \\ 0 & \beta_v \end{pmatrix} \in \mathrm{GL}_2(\mathbb{C})$.

Then given $v \notin S$, the local Langlands' $L$-functions are

$$L(s, \pi_v, r_3^0) = (1 - \alpha_v^2 \beta_v^{-1} q_v^{-s})^{-1}(1 - \alpha_v q_v^{-s})^{-1}(1 - \beta_v q_v^{-s})^{-1}(1 - \alpha_v^{-1} \beta_v^2 q_v^{-s})^{-1}$$

and

$$\begin{aligned} L(s, \pi_v, r_3) &= \det(I - r_3(t_v) q_v^{-s})^{-1} \\ &= (1 - \alpha_v^3 q_v^{-s})^{-1}(1 - \alpha_v^2 \beta_v q_v^{-s})^{-1}(1 - \alpha_v \beta_v^2 q_v^{-s})^{-1}(1 - \beta_v^3 q_v^{-s})^{-1}. \end{aligned}$$

THEOREM 2.1 ([La1], [Sh2, 3, 5]).
(1) *Fix a nontrivial additive character $\psi = \otimes_v \psi_v$ of $F\backslash\mathbb{A}$, unramified for $v \notin S$. Then at every place $v \in S$, a local $L$-function $L(s, \pi_v, r_3)$ and a local root number $\epsilon(s, \pi_v, r_3, \psi_v)$ (which equals 1 if $\pi_v$ and $\psi_v$ are unramified) can be defined so that if $L(s, \pi, r_3) = \prod_v L(s, \pi_v, r_3)$ and $\epsilon(s, \pi, r_3) = \prod_v \epsilon(s, \pi_v, r_3, \psi_v)$, then*

$$L(s, \pi, r_3) = \epsilon(s, \pi, r_3) L(1 - s, \tilde{\pi}, r_3),$$

*where $\tilde{\pi}$ is the contragredient of $\pi$. Similarly for $L(s, \pi, r_3^0)$.*

(2) $L(s, \pi, r_3)$ and $L(s, \pi, r_3^0)$ both converge absolutely for $\operatorname{Re}(s) > 1$ and are therefore holomorphic and nonzero over that region.

(3) They can be extended to meromorphic functions of $s$ on $\mathbb{C}$.

(4) They never vanish on the line $\operatorname{Re}(s) = 1$.

PROPOSITION 2.2 (Shahidi [Sh5]).

(1) Suppose $\omega_\pi = 1$. Then $L(s, \pi, r_3^0)$ and $L(s, \pi, r_3)$ are holomorphic at $s = \frac{1}{2}$.

(2) The L-function $L(s, \pi, r_3^0)$ has no poles except possibly for simple ones on the interval $\frac{1}{2} \leq s < 1$. If $\Pi$ is the Gelbart-Jacquet lift of $\pi$, then $L(s, \pi, r_3^0) = L(s, \pi \times \Pi)/L(s, \pi)$.

PROPOSITION 2.3. *Suppose $\pi$ is a monomial cuspidal representation of* $\mathrm{GL}_2(\mathbb{A})$, *i.e.*, $\pi \otimes \eta \simeq \pi$ *for some nontrivial character $\eta$ of $F^* \backslash \mathbb{A}^*$; then $L(s, \pi, r_3^0)$ has poles possibly only at $s = 0, 1$ and nowhere else. It has a pole at $s = 1$ and $0$ if and only if $\pi$ corresponds to the two-dimensional Galois representations of the Weil group of $F$ with the image isomorphic to the symmetric group $S_3$. The poles are simple.*

*Proof.* We give a proof from [Sh5], [Za]. Suppose $\pi \otimes \eta \simeq \pi$ for a nontrivial grössencharacter $\eta$. Then $\eta^2 = 1$ and $\eta$ determines a quadratic extension $E/F$. According to [L-La], there is a grössencharacter $\chi$ of $E$ such that $\pi = \pi(\chi)$, where $\pi(\chi)$ is the automorphic representation whose local factor at $v$ is the one attached to the representation of the local Weil group induced from $\chi_v$. Let $\chi'$ be the conjugate of $\chi$ by the action of the nontrivial element of the Galois group. Then the Gelbart-Jacquet lift is given by

$$\Pi = \operatorname{Ind}_{P(\mathbb{A})}^{\mathrm{GL}_3(\mathbb{A})}(\pi(\chi\chi'^{-1}) \otimes \eta),$$

where $P$ is the standard maximal parabolic subgroup of $\mathrm{GL}_3$ of type $(2,1)$. The representation $\pi(\chi\chi'^{-1})$ is cuspidal unless $\chi\chi'^{-1}$ factors through the norm. Then

$$\begin{aligned}L(s, \pi, r_3^0) &= L(s, \pi(\chi) \times \pi(\chi\chi'^{-1})) \\ &= L(s, \chi^2 \chi'^{-1}) L(s, \chi).\end{aligned}$$

This has a pole for $\operatorname{Re} s > 0$ if and only if $\pi(\chi\chi'^{-1})$ is cuspidal and $\pi(\chi\chi'^{-1}) \simeq \pi(\chi)$. This is equivalent to $\chi^3 = 1$, $\chi \neq 1$ and $\chi|_{\mathbb{A}^\times} = 1$. Representations in this class correspond to the two-dimensional Galois representations of the Weil group of $F$ whose images are isomorphic to the symmetric group $S_3$. □

*Remark* 2.1. Note that the central character of $\pi(\chi)$ is $\eta \cdot \chi|_{\mathbb{A}^\times}$. Therefore in the above case, the central character of $\pi = \pi(\chi)$ is $\omega_\pi = \eta$ and it is a nontrivial quadratic character.

So from this section on, we may assume that $\pi$ is a nonmonomial cuspidal representation of $\mathrm{GL}_2(\mathbb{A})$. Our goal is to prove that $L(s, \pi, r_3)$ and $L(s, \pi, r_3^0)$ are both entire.

PROPOSITION 2.4 ([J-S], [Sh5]). *Suppose $\pi$ is nonmonomial. Then $L(s, \pi, r_3)$ has no poles on the line $\mathrm{Re}(s) = 1$.*

Therefore it is enough to show that $L(s, \pi, r_3^0)$ is holomorphic for the region $\frac{1}{2} \leq s < 1$. Our proof will not apply to $s = 1/2$ and we therefore treat that separately.

## 3. Analysis of local intertwining operators

In this section we prove two useful propositions. We normalize the local intertwining operators by means of root numbers and $L$-functions as in [Sh1] (cf. Section 2 here). Notation will remain as in [Sh1]. In both propositions $\pi_v$ is a local component of $\pi$, although the results and the proofs are valid for any local unitary representation.

PROPOSITION 3.1. *The normalized intertwining operator $N(s, \pi_v, w)$ is holomorphic for $\mathrm{Re}\, s > 0$.*

*Proof.* For $\pi_v$ tempered and $\mathrm{Re}(s) > 0$, the local factors $L(s, \widetilde{\pi}_v, r_i)$ and $M(s, \pi_v, w)$ are holomorphic. Moreover $L(s, \widetilde{\pi}_v, r_i)$ is never zero. Next assume $\pi_v$ is a complementary series; $\pi_v = \pi(\mu|\ |^r, \mu|\ |^{-r})$, $0 < r < \frac{1}{2}$. Here we have suppressed the dependence of $\mu$ and $r$ on $v$. We follow [K2] and [Za]. Under the identification $M \simeq \mathrm{GL}_2$,

$$(3.1) \qquad \mathrm{Ind}_{P(F_v)}^{G(F_v)} \pi_v \otimes \exp(\langle s\widetilde{\alpha}, H_P(\ )\rangle) = \mathrm{Ind}_{B(F_v)}^{G(F_v)} \chi \otimes \exp(\langle \Lambda, H_B(\ )\rangle),$$

where $\Lambda = (2r)\beta_3 + (s - 3r)\beta_4$. Here $\chi$ is defined by conditions $\chi \circ \beta_3^\vee = \chi \circ \beta_5^\vee = \mu$. Then $\chi \circ \beta_1^\vee = 1$ and $\chi \circ \beta_6^\vee = \mu$. We can identify the operator $M(s, \pi_v, w)$ with $M(\Lambda, \chi, w)$.

For non-archimedean places, by Winarsky's result [W, p. 952], $M(\Lambda, \chi, w)$ is holomorphic unless $\mu = 1$, $s \in \{\pm r, \pm 3r\}$ or $\mu^2 = 1$, $s = 0$. But these are just poles of $L(s, \widetilde{\pi}_v, r_3^0)L(2s, \omega_{\pi_v}^{-1})$ and the assertion is proved. Note that $L(s, \widetilde{\pi}_v, r_3^0)L(2s, \omega_{\pi_v}^{-1})$ has no zeros. For archimedean places, we use [Sh6, p. 110]. □

The following proposition is due to [Za]. In [K2], only the case $\mathrm{Re}(s - 3r) > 0$ was treated.

PROPOSITION 3.2. *Suppose $1/2 \leq \mathrm{Re}(s) < 1$. Then for each $v$, the image of $N(s, \pi_v, w)$ is irreducible. We denote this image by $J(s, \pi_v)$. If $\pi_v$ is tempered, it is the ordinary Langlands' quotient of $I(s, \pi_v)$.*

*Proof.* Note that $-1 < \text{Re}(s - 3r) < 1$, $\langle \Lambda, \beta_6^\vee \rangle = s - 3r$, and $\langle \Lambda, \beta_1^\vee \rangle = 2r$. Therefore, the normalized operators $N(\rho_6, \Lambda)$ and $N(\rho_6, \rho_6\Lambda)$ are both holomorphic and $N(\rho_6, \rho_6\Lambda)N(\rho_6, \Lambda) = I$. (See the proof of Proposition 3.1 and the form of the normalizing factors.) Note that we have suppressed the dependence on $\chi$ in the notation. Therefore $N(\rho_6, \Lambda)$ is an isomorphism. In the same way, we can show that $N(\rho_1, \Lambda)$ is an isomorphism.

If $\text{Re}(s-3r) > 0$, then $\Lambda$ is in the positive Weyl chamber for $\frac{1}{2} \leq \text{Re}(s) < 1$ and so it is irreducible by Langlands' classification. Let us be more precise. The element $w = \rho_6\rho_1\rho_6\rho_1\rho_6$ is not the longest element of the Weyl group associated to the Borel subgroup. However, $w\rho_1 = \rho_1 w$ is. Then $N(\rho_1 w, \Lambda) = N(\rho_1, w\Lambda)N(w, \Lambda)$. By Langlands' classification, the image of $N(\rho_1 w, \Lambda)$ is irreducible and by the same method $N(\rho_1, w\Lambda)$ is an isomorphism. Therefore, the image of $N(w, \Lambda)$ is irreducible.

If $\text{Re}(s-3r) < 0$, then $\Lambda$ is not in the positive Weyl chamber. But $\rho_6(\Lambda) = (s-r)\beta_3 + (3r-s)\beta_4$ is in the positive Weyl chamber. Now use the fact that $\rho_1 w \rho_6 = \rho_6 \rho_1 w$. So $N(\rho_1 w, \rho_6\Lambda)N(\rho_6, \Lambda) = N(\rho_6, \rho_1 w\Lambda)N(\rho_1, w\Lambda)N(w, \Lambda)$. We note that the normalized operators attached to $\rho_1$ and $\rho_6$ are isomorphisms. Therefore, again the image of $N(w, \Lambda)$ is irreducible.

Suppose $\text{Re}(s - 3r) = 0$. Then by inducing in stages, we have:

$$\text{Ind}_{B(F_v)}^{G(F_v)} \chi \otimes \exp(<\Lambda, H_B(\ )>) = \text{Ind}_{P'(F_v)}^{G(F_v)} \pi'_v \otimes \exp(\langle 2r\tilde{\alpha}', H_{P'}(\ )\rangle),$$

where $\pi'_v = \text{Ind}_{B_0(F_v)}^{\text{GL}_2(F_v)} \chi(1, \mu)$ and $P'$ is the other maximal parabolic subgroup of $G$. We note from [K2] that $P' = M'N'$, $M' \simeq \text{GL}_2$ and $\tilde{\alpha}' = \beta_3$ and $\chi(1, \mu)$ is the character of $T(F)$ with respect to this identification. Here $B_0 = M' \cap B$. Again by Langlands' classification, the image of $N(s, \pi_v, w)$ is irreducible.

The fact that $J(s, \pi_v)$ is the Langlands quotient when $\pi_v$ is tempered follows from the fact that $M(s, \pi_v, w)$ and $N(s, \pi_v, w)$ are proportional if $\pi_v$ is tempered (cf. proof of Proposition 3.1). $\square$

## 4. Proof of the main theorem

Recall the following from [K3]:

OBSERVATION 4.1. *If $L(s, \pi, r_3^0)$ has a pole for $\frac{1}{2} \leq s < 1$, then $J(s, \tilde{\pi}) = \otimes_v J(s, \tilde{\pi}_v)$ belongs to the residual spectrum $L_{\text{dis}}^2(G(F)\backslash G(\mathbb{A}))_{(M,\tilde{\pi})}$, and consequently each $J(s, \tilde{\pi}_v)$ is unitary.*

*Proof.* By [La2], [M-W1] the global intertwining operator $M(s, \pi, w)$ gives the poles of the Eisenstein series. In view of the results in Section 2 and Proposition 3.1, they are precisely those of $L(s, \pi, r_3^o)$ if $1/2 \leq s < 1$. Moreover the residue at a given (simple) pole $s$ on this interval is $J(s, \tilde{\pi}) = \otimes_v J(s, \tilde{\pi}_v)$, an irreducible representation (Proposition 3.2) which by [La2], [M-W1] appears in $L^2$ and is therefore unitary. More precisely by [La1], [M-W1], the kernel of the residue of the Eisenstein series, a map from $I(s, \tilde{\pi})$ into the $L_2$-space, is exactly that of the normalized operator $N(s, \tilde{\pi}, w) = \otimes_v N(s, \tilde{\pi}_v, w)$ whose image is the irreducible representation $J(s, \tilde{\pi})$. □

The new local result which we need to rely on is in a recent article of Muić [Mu] in which he has determined the unitary dual of $G_2$ over a non-archimedean field completely. Here, we recall a part of the classification needed.

THEOREM 4.2 (Muić [Mu]). *If $\pi_v$ is tempered and $s$ is real, then $J(s, \pi_v)$ is unitary if and only if one of the following conditions is satisfied:*

(1) $\pi_v$ is a supercuspidal representation $\rho$ such that $\rho \simeq \tilde{\rho}$, $\omega_\rho = 1$, and $0 < s \leq \frac{1}{2}$.
(2) $\pi_v$ is a supercuspidal representation $\rho$ such that $\rho \simeq \tilde{\rho}$, and $\text{Im}(\tau) \simeq S_3$ (the symmetric group in 3 variables), and $0 < s \leq 1$, where $\rho = \pi(\tau)$ and $\tau: W_F \longmapsto GL(2, \mathbb{C})$ is the attached admissible homomorphism. Then $\det \tau = \omega_\rho$ via class field theory.
(3) $\pi_v$ is a non-supercuspidal discrete series and $0 < s \leq \frac{1}{2}$.
(4) $\pi_v \simeq \pi(\mu_v, \mu_v^{-1})$, $\mu_v^3 \neq 1$, and $0 < s \leq \frac{1}{2}$.
(5) $\pi_v \simeq \pi(\mu_v, \mu_v^{-1})$, $\mu_v^3 = 1$, and $0 < s \leq \frac{1}{2}$ or $s = 1$.
(6) $\pi_v \simeq \pi(1, \mu_v)$, $\mu_v$ has order two, and $0 < s \leq 1$.

Next we observe that if $\pi_v = \pi(|\ |^r, |\ |^{-r})$, $0 < r < 1/2$, is a complementary series representation of $M(F_v) = \text{GL}_2(F_v)$, then by (3.1):

$$I(s, \pi_v) = I(\Lambda, 1),$$

where $\Lambda = 2r\beta_3 + (s - 3r)\beta_4$ is as in Proposition 3.1. By Theorem 5.2 of [Mu], the shaded upper triangle of the unitary dual coming from the conjugacy class of Borel subgroups is defined by inequalities $\langle \Lambda, \beta_2^\vee \rangle > 1$, $\langle \Lambda, \beta_3^\vee \rangle < 1$, $\langle \Lambda, \beta_1^\vee \rangle > 0$, and $\langle \Lambda, \beta_6^\vee \rangle > 0$. Since $\langle \Lambda, \beta_2^\vee \rangle = s + 3r$ and $\langle \Lambda, \beta_3^\vee \rangle = s + r$, while $\langle \Lambda, \beta_1^\vee \rangle = 2r$ and $\langle \Lambda, \beta_6^\vee \rangle = s - 3r$, this triangle has the points $(1/6, 1/2)$, $(1/4, 3/4)$, and $(0, 1)$ as its vertices in the $rs$-plane. Consequently for each $0 < r < 1/4$, there exists $1/2 < s < 1$ for which $J(s, \pi_v)$ is unitary and all such $s$ are attained for $0 < r < 1/4$.

In view of this observation it is clear that to eliminate the poles between $1/2$ and $1$ (Theorem 4.9), beside those of Section 2 and Observation 4.1, one still needs another global ingredient since *a priori* local components $\pi_v$ could all be in the complementary series, i.e. nontempered. For this we use a recent result of Ramakrishnan [Ra] that allows an arbitrary cuspidal representation to have tempered components at a subset of finite places of $F$ with a fairly large lower Dirichlet density, a very deep result. In both proofs that we give for our main result (Theorem 4.9), we need his result, but to differing degrees. One requires his density, but not fully, while the other uses the fact that the cusp form has at least one nondiscrete (non-supercuspidal is enough) tempered component, still a fairly deep global result. Consequently the second proof relies more on local results.

We recall Ramakrishnan's result [Ra]:

THEOREM 4.3 (Ramakrishnan [Ra]). *Let $\pi$ be a cuspidal, unitary representation of* $\mathrm{GL}_2(\mathbb{A})$. *Let $S(\pi)$ be the set of primes where $\pi_v$ is tempered. For any set $X$ of primes, denote by $\underline{\delta}(X)$ the lower Dirichlet density of $X$. Then*

$$\underline{\delta}(S(\pi)) \geq \frac{9}{10}.$$

Next, let $J(\nu, \sigma)$ be the Langlands quotient for an irreducible unitary representation of $G$. Let $P_1 = M_1 N_1$ be the corresponding standard parabolic subgroup, $N_1 \subset U$, so that $\sigma$ is a tempered representation of $M_1(F_v)$ and $\nu \in \mathfrak{a}_{M_1}^*$, where $\mathfrak{a}_{M_1}$ is the real Lie algebra of the split component $A_1$ of $M_1$. Then since $J(\nu, \sigma)$ is unitary, it is hermitian [K3]. This means that there exists an element $w_1$ in the Weyl group of $A_1$ in $G$ satisfying $w_1(\sigma) \cong \sigma$ and $w_1(\nu) = -\nu$, where $\nu$ is real. We record this observation from [K3] as:

OBSERVATION 4.4. *Any unitary representation is hermitian.*

LEMMA 4.5. *Assume $\pi_v$ is tempered and $J(s, \pi_v)$ is unitary. Moreover suppose $s$ is real. Then $\pi_v \cong \tilde{\pi}_v$. In particular $\omega_{\pi_v}^2 = 1$.*

*Proof.* Since $\pi_v$ is tempered and $J(s, \pi_v)$ is unitary, it must be hermitian by Observation 4.4. Consequently $w(\pi_v) \cong \pi_v$. By the discussion in the last part of the proof of Proposition 6.1 of [Sh5], $w(\pi_v) = \pi_v \otimes \omega_{\pi_v}^{-1} = \tilde{\pi}_v \cong \pi_v$. □

LEMMA 4.6. *Suppose $L(s, \pi, r_3^0)$ has a pole for $1/2 \leq s < 1$. Then $\omega_\pi^2 = 1$.*

*Proof.* By Ramakrishnan's Theorem 4.3 and Lemma 4.5, $\omega_{\pi_v}^2 = 1$ for a set of primes whose lower Dirichlet density is at least $9/10$. One can now apply either the equidistribution of values of $\omega_\pi$, or as pointed out to us by Ramakrishnan, use the simpler result of Hecke that two idele class characters agreeing at all the places in a set of Dirichlet density larger than $1/2$ are equal, to conclude that $\omega_\pi^2 = 1$. □

LEMMA 4.7. *The L-function $L(s, \pi, r_3^0)$ is holomorphic at $s = 1/2$.*

*Proof.* By Part 1 of Proposition 2.2, we may assume $\omega_\pi \neq 1$. By Lemma 4.6, $\omega_\pi^2 = 1$, $\omega_\pi \neq 1$. Let $E/F$ be the quadratic extension of $F$ defined by $\omega_\pi$.

Denote by $\Pi$ the base change [La3] of $\pi$ to $\mathrm{GL}_2(\mathbb{A}_E)$. Since $\pi$ is not monomial $\Pi$ remain cuspidal [La3]. Then

$$L(s, \Pi, r_3^0) = L(s, \pi, r_3^0) L(s, \pi \otimes \omega_\pi, r_3^0).$$

Note that $\pi \otimes \omega_\pi = \tilde{\pi}$ and by the functional equation (Theorem 2.1.1),

$$L(1/2, \pi, r_3^0) = \varepsilon(1/2, \pi, r_3^0) L(1/2, \tilde{\pi}, r_3^0).$$

Thus if $L(s, \pi, r_3^0)$ has a pole at $s = 1/2$, then $L(s, \Pi, r_3^0)$ must have a double pole at $s = 1/2$. This contradicts the holomorphy (in fact simplicity of the pole is enough) of $L(s, \Pi, r_3^0)$ at $s = 1/2$, as the central character of $\Pi$ is trivial, again by Proposition 2.2.1, which completes the lemma. □

The following lemma relies more on local results of Muić [Mu], and Propositions 3.1 and 3.2 here. It reduces our use of [Ra] to $\pi$ having a single nondiscrete tempered component, still a fairly deep global result. We will use this lemma to give a second proof of our main result.

LEMMA 4.8. *$J(s, \pi_v)$ is hermitian if and only if $\pi_v \cong \tilde{\pi}_v$. Consequently if $L(s, \pi, r_3^0)$ has a pole for $s > 0$, then $\pi \cong \tilde{\pi}$.*

*Proof.* In view of Lemma 4.5, we may assume $\pi_v = \pi(\mu|\ |^r, \mu|\ |^{-r})$, where $\mu$ is a unitary character of $F_v^*$ and $0 < r < 1/2$. Again we have suppressed the dependence of $\mu$ and $r$ on $v$. Then

$$I(s, \pi_v) \cong \mathrm{Ind}_{B(F_v)}^{G(F_v)} \chi \otimes \exp(\langle \Lambda, H_B(\ )\rangle)$$
$$\stackrel{\mathrm{def}^n}{=} I(\Lambda, \chi),$$

where $\Lambda = (2r)\beta_3 + (s - 3r)\beta_4$ and $\chi$ is as in the proof of Proposition 3.1; i.e. $\chi \circ \beta_3^\vee = \chi \circ \beta_5^\vee = \mu$. The longest element in the Weyl group of $T$ in $G$ sends $\chi$ to $\chi^{-1}$. Suppose $s - 3r > 0$. Then the Langlands quotient $J(s, \pi_v)$ of (3.1) being hermitian implies $\chi^{-1} = \chi$, leading to $\omega_{\pi_v} = \mu^2 = 1$. Now suppose $\mathrm{Re}(s - 3r) < 0$. Then as in Proposition 3.2, the Langlands quotient of the (standard) module induced from $\chi' = \rho_6(\chi)$ and $\rho_6(\Lambda)$ is again hermitian. This time $\chi' \circ \beta_5^\vee = 1$, while $\chi' \circ \beta_1^\vee = \mu$. The same argument now applies. The case $s - 3r = 0$ is even easier. The lemma is now proved. □

We are now ready to prove the main result of our paper:

THEOREM 4.9. *Suppose $\pi$ is not monomial. Then the symmetric cube L-function $L(s, \pi, r_3)$ and the adjoint cube L-function $L(s, \pi, r_3^0)$ are both entire.*

*Proof.* We shall give three proofs. The first one relies more on global results [Ra] than local ones [Mu]. Besides the results in Section 2, it requires Lemmas 4.5, 4.6, and 4.7. The second one is more local and uses more of classification of the unitary dual of $G_2$ over a non-archimedean field [Mu]. But it still needs Ramakrishnan's global result [Ra], although not as strongly. Besides the results in Sections 2 and 3, it requires Lemma 4.8. □

*First proof.* In view of the results in Section 2 and Lemma 4.7 we only need to show that the $L$-functions are holomorphic for $1/2 < s < 1$. Suppose $L(s, \pi, r_3^0)$ has a pole for $1/2 < s < 1$. Then by Lemma 4.6, $\omega_\pi^2 = 1$. By equidistribution of values of $\omega_\pi$, $\omega_{\pi_v} = 1$ for a set of (finite) primes $v$ with a Dirichlet density at least $1/2$. On the other hand by Theorem 4.3, the lower Dirichlet density for primes with $\pi_v$ tempered is at least $9/10$. By Observation 4.1 we therefore have a finite place $v$ and a nondiscrete tempered component $\pi_v$ with a trivial central character such that $J(s, \pi_v)$ is unitary. By parts (4) and (5) of Theorem 4.2 this is a contradiction, which completes the theorem. □

*Second proof.* Assume that $s$ (with $1/2 < s < 1$) is a pole of $L(s, \pi, r_3^0)$. Then by Lemma 4.8, $\pi \cong \tilde\pi = \pi \otimes \omega_\pi^{-1}$. Unless $\omega_\pi = 1$, $\pi$ is monomial. We therefore need to assume $\omega_\pi = 1$. By Theorem 4.3 there is at least one unramified tempered component $\pi_v$. Since $\omega_{\pi_v} = 1$, Parts 4 and 5 of Theorem 4.2 show that $J(s, \pi_v)$ is not unitary, contradicting the pole at $s$. The pole at $s = 1/2$ is eliminated by Proposition 2.2.1 and the theorem is now a consequence of the results of Section 2, particularly the functional equation.□

*Sketch of a third proof.* Our global ingredient is now the much weaker result that given $\varepsilon > 0$, there exists a local component $\pi_v = \pi(\mu_v| \ |_v^{r_v}, \mu_v| \ |_v^{-r_v})$, $0 \leq r_v < 1/2$, $\mu_v \in \hat F_v^*$, such that $r_v < \varepsilon$ (cf. [K3]). We now drop the subscript $v$ from our notation and observe that if $\chi_1$ and $\chi_2$ are characters of $F_v^*$ as in Theorem 5.2 of [Mu], then $\chi_1 = \chi \cdot \beta_5^\vee = \mu = \mu_v$ and $\chi_2 = \chi . \beta_1^\vee = 1$, where $\chi$ is as in the proof of Proposition 3.1. We translate the unitary dual of $G_2$ coming from the conjugacy class of its Borel subgroups, given by Theorem 5.2 of [Mu], to our $rs$-plane, $r = r_v$, as in the discussion immediately after Theorem 4.2. If $\mu = 1$, then $\chi_1 = \chi_2 = 1$ and the contribution to the unitary dual will consist of the two shaded triangles of Figure 1, page 483, of [Mu], which in the $rs$-plane have $(0,1)$, $(1/4, 3/4), (1/6, 1/2)$, and $(0, 1/2), (1/6, 1/2), (0,0)$ as their vertices. On the other hand, if $\mu^2 = 1$, but $\mu \neq 1$, then $\chi_1^2 = \chi_2 = 1$, but $\chi_1 \neq 1$, and we are in Case (iii) of Theorem 5.2 of [Mu]. The contribution now comes only from the lower triangle $(0, 1/2), (1/6, 1/2), (0, 0)$. Now, if $s$ is a pole for $L(s, \pi, r_3^0), 1/2 < s < 1$, we choose $v$ such that $0 < r_v < (1-s)/3$. This gives a parameter inside the triangle whose vertices are $(0,1)$, $(1/6, 1/2)$, and $(0,1/2)$ which cannot be unitary, a contradiction. Observe that the fact that the inside of this last triangle cannot afford any unitary parameter is crucial and remarkable. To conclude we apply Lemma 4.7. □

*Remark* 4.10. As discussed in the introduction, using an integral representation, Bump, Ginzburg, and Hoffstein [B-G-H], have also proved the holomorphy of the partial $L$-function $L_S(s, \sigma, r_3)$ for $\text{Re}(s) > 3/4$, provided that $F$ has cube roots of unity.

*Remark* 4.11. From Proposition 2.2 and Theorem 4.9, it follows that $L(s, \pi \times \Pi)/L(s, \pi)$ is entire, where $\Pi$ is the Gelbart-Jacquet lift of $\pi$. This means that the zeros of $L(s, \pi)$ are among those of $L(s, \pi \times \Pi)$.

*Remark* 4.12. From Remark 1.1 and Theorem 4.9, it follows that the partial Rankin triple $L$-function $L_S(s, \pi \times \pi \times \pi)$ is entire for a nonmonomial representation $\pi$. It also follows that the zeros of $L_S(s, \pi \otimes \omega_\pi)^2$ are among those of $L_S(s, \pi \times \pi \times \pi)$ and in particular $L_S(s, \pi \times \pi \times \pi)$ could have double zeros at $s = 1/2$.

*Remark* 4.13.
(1) In calculation of the residual spectrum in [K2], the first author had to assume the holomorphy of $L(s, \pi, r_3^0)$ for nonmonomial representations, together with a bound on the Fourier coefficients, in order to prove Proposition 3.2. Since we now have these facts by Proposition 3.2 and Theorem 4.9, we can remove those assumptions from [K2].
(2) In calculation of the residual spectrum associated to Borel subgroups, in [K2] the first author assumed that the archimedean components were spherical. Now Zampera [Za] has removed this condition and therefore the archimedean places can be treated the same way as others.
(3) The first author likes to make the following corrections to [K2]. On page 1245, line 10, $G$ should be assumed to have anisotropic center. In line 12 of page 1248, Zampera should be referred to with a masculine pronoun. In the character table of $S_4$ on page 1266, $\psi_{211}$, $\psi_{22}$, and $\psi_{31}$, should respectively be changed to $\psi_{22}$, $\psi_{31}$, and $\psi_{211}$. Finally, in line 1 of page 1271, $\Lambda$ should be $\Lambda = (s - r)\beta_3 + (2r)\beta_4$ and therefore the assumption $r < 1/6$ in Theorem 5.1 is unnecessary.

SOUTHERN ILLINOIS UNIVERSITY, CARBONDALE, IL
*E-mail address*: henrykim@math.siu.edu

PURDUE UNIVERSITY, WEST LAFAYETTE, IN
*E-mail address*: shahidi@math.purdue.edu

REFERENCES

[Bo]    A. BOREL, Automorphic $L$-functions, in *Automorphic Forms, Representations and L-Functions*, Proc. Sympos. Pure Math. **33**; Part II, Amer. Math. Soc., Providence, RI, 1979, 27–61.

[B-D-H-I]    D. BUMP, W. DUKE, J. HOFFSTEIN, and H. IWANIEC, An estimate for the Hecke eigenvalues of Maass forms, Internat. Math. Res. Notices **4** (1992), 75–81.

[B-G-H]    D. BUMP, D. GINZBURG, and J. HOFFSTEIN, Symmetric cube $L$-functions, Invent. Math., to appear.

[C-S]    W. CASSELMAN and J. A. SHALIKA, The unramified principal series of $p$-adic groups II, The Whittaker function, Compositio Math. **41** (1980), 207–231.

[Co-PS]    J. W. COGDELL and I. I. PIATETSKI-SHAPIRO, Converse theorems for $GL_n$, Publ. Math. I.H.E.S. **79** (1994), 157–214.

[D-I]    W. DUKE and H. IWANIEC, Estimates for coefficients of $L$-functions IV, Amer. J. Math. **116** (1994), 207–217.

[De]    P. DELIGNE, Valeurs de functions $L$ et périodes d'intégrales, in *Automorphic Forms, Representations and L-Functions*, Proc. Sympos. Pure Math. **33**; Part II, Amer. Math. Soc., Providence, RI, 1979, 313–346.

[F-Gol]    S. FRIEDBERG AND D. GOLDBERG, On local coefficients for nongeneric representations of some classical groups, Compositio Math. **116** (1999), 133–166.

[Ga]    P. GARRETT, Decomposition of Eisenstein series: Rankin triple product, Ann. of Math. **125** (1987), 209–237.

[Ge-J]    S. GELBART AND H. JACQUET, A relation between automorphic representations of GL(2) and GL(3), Ann. Sci. École Norm. Sup. **11** (1978), 471–552.

[Go]    R. GODEMENT, Formes automorphes et produits eulerian (d'apres R. P. Langlands), Seminaire Bourbaki **349** (1968–69), 38–53.

[H-Ra]    J. HOFFSTEIN and D. RAMAKRISHNAN, Siegel zeros and cusp forms, Internat. Math. Res. Notices **6** (1995), 279–308.

[Ik]    T. IKEDA, On the location of poles of the triple $L$-functions, Compositio Math. **83** (1992), 187–237.

[J-S]    H. JACQUET and J. A. SHALIKA, On Euler products and the classification of automorphic forms II, Amer. J. Math. **103** (1981), 777–815.

[K1]    H. KIM, The residual spectrum of $Sp_4$, Compositio Math. **99** (1995), 129–151.

[K2]    ———, The residual spectrum of $G_2$, Canad. J. Math. **48** (1996), 1245–1272.

[K3]    ———, Langlands-Shahidi method and poles of automorphic $L$-functions: Application to exterior square $L$-functions, Canad. J. Math., to appear.

[L-La]    J-P. LABESSE and R. P. LANGLANDS, $L$-indistinguishability for SL(2), Canad. J. Math. **31** (1979), 726–785.

[La1]    R. P. LANGLANDS, *Euler Products*, Yale University Press, New Haven, Conn., 1971.

[La2]    ———, *On the Functional Equations Satisfied by Eisenstein Series*, Lecture Notes in Math. **544**, Springer-Verlag, New York, 1976.

[La3]    ———, *Base Change for* GL(2), Ann. of Math. Studies **96**, Princeton University Press, Princeton, NJ, 1980.

[La4]    ———, Problems in the theory of automorphic forms, Lecture Notes in Math. **170**, Springer-Verlag, New York, 1970, 18–86.

[Lu-R-Sa1]    W. LUO, Z. RUDNICK, and P. SARNAK, On Selberg's eigenvalue conjecture, Geom. Funct. Anal. **5** (1995), 387–401.

[Lu-R-Sa2]    ———, On the generalized Ramanujan conjecture for GL($n$), preprint, 1996.

[M-W1]    C. MŒGLIN and J.-L. WALDSPURGER, *Spectral Decomposition and Eisenstein Series. Une Paraphrase de l'Ecriture*, Cambridge Tracts in Math. **113**, Cambridge University Press, Cambridge, 1995.

[M-W2]    ———, Le spectre résiduel de GL($n$), Ann. Sci. Ećole Norm. Sup. **22** (1989), 605–674.

[Mo-Sh]    C. J. MORENO AND F. SHAHIDI, The $L$-function $L_3(s, \pi_\Delta)$ is entire, Invent. Math. **79** (1985), 247–251.

[Mu]    G. MUIĆ, The unitary dual of $p$-adic $G_2$, Duke Math. J. **90** (1997), 465–493.

[Muk]   V. K. Murty, On the Sato-Tate conjecture, in *Number Theory Related to Fermat's Last Theorem*, Birkhauser-Verlag, Boston, MA, 1982, 195–205.

[Mur]   R. Murty, Oscillations of Fourier coefficients of modular forms, Math. Ann. **262** (1983), 431–446.

[PS-Ral]   I. I. Piatetski-Shapiro and S. Rallis, Rankin triple $L$-functions, Compositio Math. **64** (1987), 31–115.

[Ra]   D. Ramakrishnan, On the coefficients of cusp forms, Math. Res. Lett. **4** (1997), 295–307.

[Se]   J-P. Serre, *Abelian $\ell$-Adic Representations and Elliptic Curves*, W. A. Benjamin, Inc, New York, 1968; second edition, Addison-Wesley, Redwood City, CA, 1989.

[Sh1]   F. Shahidi, A proof of Langlands' conjecture on Plancherel measures; Complementary series for $p$-adic groups, Ann. of Math. **132** (1990), 273–330.

[Sh2]   ———, Functional equation satisfied by certain $L$-functions, Compositio Math. **37** (1978), 171–207.

[Sh3]   ———, On certain $L$-functions, Amer. J. Math. **103** (1981), 297–355.

[Sh4]   ———, On the Ramanujan conjecture and finiteness of poles for certain $L$-functions, Ann. of Math. **127** (1988), 547–584.

[Sh5]   ———, Third symmetric power $L$-functions for GL(2), Compositio Math. **70** (1989), 245–273.

[Sh6]   ———, Whittaker models for real groups, Duke Math. J. **47** (1980), 99–125.

[Sh7]   ———, Symmetric power $L$-functions for GL(2), in *Elliptic Curves and Related Topics* (H. Kisilevsky and M. R. Murty, eds.), CRM Proc. Lectures Notes **4**, Amer. Math. Soc. Providence, RI, 1994, 159–182.

[S-Tan]   J. A. Shalika and S. Tanaka, On an explicit construction of a certain class of automorphic forms, Amer. J. Math. **91** (1969), 1049–1076.

[Shi]   G. Shimura, On the holomorphy of certain Dirichlet series, Proc. London Math. Soc. **31** (1975), 79–98.

[Sp]   B. Speh, Unitary representations of $GL(n, \mathbb{R})$ with nontrivial $(\mathfrak{g}, K)$-cohomology, Invent. Math. **71** (1983), 443–465.

[Ta]   M. Tadić, Classification of unitary representations in irreducible representations of general linear groups (non-archimedean case), Ann. Sci. École Norm. Sup. **19** (1986), 335–382.

[T]   J. T. Tate, Algebraic cycles and poles of zeta functions, *Arithmetical Algebraic Geometry* (Proc. Conf. Purdue Univ., 1963), Harper & Row, New York, 1965, 93–110.

[V]   D. A. Vogan, The unitary dual of $G_2$, Invent. Math. **116** (1994), 677–791.

[W]   N. Winarsky, Reducibility of principal series representations of $p$-adic Chevalley groups, Amer. J. Math. **100** (1978), 941–956.

[Za]   S. Žampera, The residual spectrum of the group of type $G_2$, J. Math. Pures Appl. **76** (1997), 805–835.